\documentclass[11pt,a4paper]{article}
\usepackage[T1]{fontenc}
\usepackage{lmodern}
\usepackage[a4paper,margin=27mm]{geometry}
\usepackage{amsmath,amssymb,amsthm,mathtools}
\usepackage{booktabs,array,longtable}
\usepackage[colorlinks=true,linkcolor=blue,citecolor=blue,urlcolor=blue]{hyperref}
\numberwithin{equation}{section}
\newtheorem{theorem}{Theorem}[section]
\newtheorem{proposition}[theorem]{Proposition}
\newtheorem{lemma}[theorem]{Lemma}
\newtheorem{corollary}[theorem]{Corollary}
\theoremstyle{remark}\newtheorem{remark}[theorem]{Remark}
\newcommand{\R}{\mathbb R}
\newcommand{\g}{\mathfrak g}
\newcommand{\vol}{\mathrm{vol}}
\newcommand{\dd}{\mathrm d}
\newcommand{\spn}{\operatorname{span}}
\newcommand{\Ann}{\operatorname{Ann}}
\newcommand{\tr}{\operatorname{tr}}
\newcommand{\ad}{\operatorname{ad}}
\newcommand{\cH}{\mathcal H}
\newcommand{\cU}{\mathcal U}
\newcommand{\cB}{\mathcal B}

\title{Left-Invariant Kernels of the Hodge--de Rham Operator\\on Three-Dimensional Lorentzian Lie Groups}
\author{FANQI LIN}\date{}
\begin{document}
\maketitle
\vspace{-2.6em}
\begin{abstract}
We determine the kernel of the Hodge--de Rham operator on left-invariant differential forms for the seven customary families of three-dimensional Lorentzian metric Lie algebras. Every parameter stratum is treated, including zero structure constants, singular bracket matrices, and null modular directions. If $\theta(X)=\tr\ad_X$ is the modular character, the full kernel is the direct sum of constants, the space $K=\Ann([\g,\g])\cap\theta^\perp$, its Hodge dual, and the volume forms. We give explicit real bases in all seven families and compute the exterior differential and geometric codifferential on every basis form. The non-unimodular calculation remains valid without the customary extra normalization conditions, and hence also covers metric representatives omitted by those conditions. The resulting dimensions are $2,4,6$, or $8$; the two middle-degree dimensions agree. Explicit examples distinguish the kernel of the first-order operator from the kernel of its square. All assertions concern smooth left-invariant forms, with no square-integrability condition.
\end{abstract}


\section{Introduction}
The spin Dirac operator and the Hodge--de Rham operator provide distinct first-order geometric equations on a metric Lie group. The former acts on a spinor module, whereas the latter acts on the complete exterior algebra. Gil-Garc\'ia and Russo \cite{GilRusso} study left-invariant harmonic spinors in dimension three, in both Riemannian and Lorentzian signature. Their problem suggests asking for an equally explicit treatment of
\begin{equation}\label{eq:Dintro}
 D=\dd+\delta:\Lambda^\bullet\g^*\longrightarrow\Lambda^\bullet\g^*
\end{equation}
on left-invariant differential forms. The passage between these two problems is not a substitution in the spinor matrix: the representation space and the induced connection both change.

There is also a basic change in the existence question. Constants and volume forms always belong to the left-invariant kernel of $D$. Thus the existence of a nonzero kernel does not restrict the metric. The substantial questions are the existence of middle-degree solutions, their explicit representatives, and the dependence of the kernel dimension on the structure constants. In Lorentzian signature one must additionally distinguish $\ker D$ from $\ker D^2$: a positive-definite Hodge energy argument is unavailable.

Our principal observation is that, in dimension three, the full inhomogeneous equation reduces to a single space of one-forms. Let
\begin{equation}\label{eq:introK}
 \theta(X)=\tr(\ad_X),\qquad
 W=\Ann([\g,\g]),\qquad K=W\cap\theta^\perp,
\end{equation}
where the orthogonal complement is taken in $\g^*$ with its inverse metric. We prove
\begin{equation}\label{eq:mainintro}
 \ker D=\R\,1\oplus K\oplus *K\oplus\R\,\vol.
\end{equation}
This formula determines actual representatives, not only dimensions. In a pseudo-orthonormal coframe of signature $(+,+,-)$, every generator $Ae^1+Be^2+Ce^3$ of $K$ produces the explicit two-form $Ae^{23}-Be^{13}-Ce^{12}$.

We apply this reduction to the seven standard bracket families $\g_1,\ldots,\g_7$ used in \cite{BatatOnda,Wang}. The first four are unimodular; the last three have, respectively, a spacelike, Lorentzian, or degenerate unimodular kernel. The computations yield a complete list for every parameter choice in these displayed families. In particular, $\g_3$ requires eight zero/nonzero strata, $\g_4$ requires four strata, and $\g_5,\g_6$ require a rank-one subcase with an explicit choice of generator.

A qualification about the standard families is essential. The non-unimodular normalizations frequently appended to the seven-family list do not all preserve every Lorentzian metric representative. Brozos-V\'azquez, Caeiro-Oliveira, and Garc\'ia-R\'io \cite[\S3.2.2, Remarks 4.15 and 4.23]{BCG} explain the problem for the Lorentzian and degenerate unimodular kernels. Our formulas for $\g_6$ and $\g_7$ do not require those extra restrictions. We give the customary parameter conditions to fix the family notation, but prove the results on the larger unrestricted bracket families as well. Consequently the kernel computation itself is independent of the disputed normalizations.

Section~\ref{sec:reduction} fixes all signs and proves the reduction theorem. Section~\ref{sec:classification} records the seven families and the classification table. Sections~\ref{sec:unimod} and \ref{sec:nonunimod} compute every case directly. Section~\ref{sec:geometry} gives invariant consequences and coverage beyond normalization. Section~\ref{sec:square} distinguishes the first-order kernel from that of the Hodge Laplacian. The appendix provides full matrix data and an exact verification procedure. The work is a kernel classification within an established metric-Lie-algebra framework, rather than a new classification of Lorentzian metrics.

\section{Conventions and the three-dimensional reduction}\label{sec:reduction}
\subsection{The Hodge star and geometric codifferential}
Let $G$ be a connected oriented three-dimensional Lie group with a left-invariant Lorentzian metric $g$. Choose an oriented pseudo-orthonormal frame with dual coframe, denoted by $(e_i)$ and $(e^i)$, respectively, so that
\begin{equation}\label{eq:metric}
 g(e_i,e_j)=\varepsilon_i\delta_{ij},\qquad
 (\varepsilon_1,\varepsilon_2,\varepsilon_3)=(1,1,-1),\qquad
 \vol=e^{123}.
\end{equation}
We abbreviate $e^{ij}=e^i\wedge e^j$. All coefficients of left-invariant forms are constant. Their differential is the Chevalley--Eilenberg differential
\begin{equation}\label{eq:CE}
 \dd e^k(e_i,e_j)=-e^k([e_i,e_j]),\qquad
 \dd(\xi\wedge\zeta)=\dd\xi\wedge\zeta+(-1)^{\deg\xi}\xi\wedge\dd\zeta.
\end{equation}
The Hodge star is defined by $\xi\wedge *\zeta=\langle\xi,\zeta\rangle\vol$, where the pairing on forms is induced by $g^{-1}$. It acts by
\begin{align}
 *1&=\vol,& *\vol&=-1,\label{eq:star0}\\
 *e^1&=e^{23},& *e^2&=-e^{13},& *e^3&=-e^{12},\label{eq:star1}\\
 *e^{12}&=e^3,& *e^{13}&=e^2,& *e^{23}&=-e^1.\label{eq:star2}
\end{align}
Thus $*^2=-I$ in every degree. The geometric codifferential is
\begin{equation}\label{eq:delta}
 \delta|_{\Lambda^k}=(-1)^{3(k+1)+2}*\dd*,\qquad
 \delta_1=*\dd*,\quad \delta_2=-*\dd*,\quad \delta_3=*\dd*.
\end{equation}
This is the formal codifferential determined by the Lorentzian metric and its volume density. Accordingly $D$ is the Hodge--de Rham operator in this convention, and $D^2=\dd\delta+\delta\dd$.

\subsection{The modular character}
The modular character $\theta\in\g^*$ is defined in \eqref{eq:introK}. It is closed because
\begin{equation}\label{eq:theta-closed}
 \theta([X,Y])=\tr[\ad_X,\ad_Y]=0.
\end{equation}
For a left-invariant vector field $X$, differentiation of the left-invariant volume form along its flow gives
\begin{equation}\label{eq:div}
 \mathcal L_X\vol=-\theta(X)\vol,\qquad \operatorname{div}X=-\theta(X).
\end{equation}
Indeed, the Lie derivative on the three covector factors contributes the negative trace of $\ad_X$. Since $\delta\xi=-\operatorname{div}(\xi^\sharp)$, we obtain:
\begin{lemma}\label{lem:modular}
For every left-invariant one-form $\xi$,
\begin{equation}\label{eq:modident}
 \delta\xi=\langle\theta,\xi\rangle,\qquad
 \dd(*\xi)=-\langle\theta,\xi\rangle\vol,\qquad
 \delta(*\xi)=*\dd\xi.
\end{equation}
Moreover, $D1=D\vol=0$.
\end{lemma}
\begin{proof}
The first formula follows from \eqref{eq:div}. For the second, use $*\xi=\iota_{\xi^\sharp}\vol$ and Cartan's identity:
$\dd\iota_{\xi^\sharp}\vol=\mathcal L_{\xi^\sharp}\vol$, since $\dd\vol=0$.
For the third, apply \eqref{eq:delta} to the two-form $*\xi$ and use $*^2\xi=-\xi$. Finally $\dd1=0$, while $*\vol=-1$ is constant, so $\delta\vol=0$.
\end{proof}

\begin{remark}\label{rem:adjoint}
The geometric $\delta$ must not be replaced by the finite-dimensional adjoint $\dd^\dagger$ on $\Lambda^\bullet\g^*$. In general,
\begin{equation}\label{eq:correction}
 \delta=\dd^\dagger+\iota_{\theta^\sharp}.
\end{equation}
For example, $\dd^\dagger:\Lambda^1\g^*\to\R$ vanishes, whereas $\delta\xi$ need not vanish. Our calculations use \eqref{eq:delta} directly, so no assumption about the adjoint on the invariant subspace enters the proofs. The appendix verifies \eqref{eq:correction} in coordinates in the present dimension.
\end{remark}

\subsection{The full inhomogeneous kernel}
Write $\cH^k_D=\ker D\cap\Lambda^k\g^*$ and $\cH_D=\ker(D|_{\Lambda^\bullet\g^*})$. These are real vector spaces. We call their elements \emph{$D$-harmonic}; the term does not assert $L^2$ integrability or identify this kernel with that of the wave-type Hodge Laplacian.
\begin{theorem}[Kernel formula]\label{thm:reduction}
For any oriented three-dimensional Lorentzian metric Lie algebra, put
\begin{equation}\label{eq:Kdefinition}
 K=\{\xi\in\g^*: \xi([\g,\g])=0,\ \langle\theta,\xi\rangle=0\}.
\end{equation}
Then
\begin{equation}\label{eq:kernel}
 \cH_D=\R\,1\oplus K\oplus *K\oplus\R\,\vol,
\end{equation}
and every element of this kernel is separately closed and coclosed. In particular,
\begin{equation}\label{eq:dimensions}
 (\dim\cH_D^0,\dim\cH_D^1,\dim\cH_D^2,\dim\cH_D^3)
 =(1,k,k,1),\qquad \dim\cH_D=2+2k,
\end{equation}
where $k=\dim K$.
\end{theorem}
\begin{proof}
Every invariant inhomogeneous form has a unique expression
\begin{equation}\label{eq:inhom}
 \omega=f+\xi+*\zeta+h\vol,\qquad f,h\in\R,\quad \xi,\zeta\in\g^*.
\end{equation}
Lemma~\ref{lem:modular} gives its image with all degrees displayed:
\begin{equation}\label{eq:Dinhom}
 D\omega=\underbrace{\langle\theta,\xi\rangle}_{0\text{-form}}
 +\underbrace{*\dd\zeta}_{1\text{-form}}
 +\underbrace{\dd\xi}_{2\text{-form}}
 -\underbrace{\langle\theta,\zeta\rangle\vol}_{3\text{-form}}.
\end{equation}
The four summands have distinct degrees and cannot cancel. The star is invertible, and $\dd\xi=0$ precisely when $\xi$ annihilates the derived algebra. Thus $D\omega=0$ is equivalent to $\xi,\zeta\in K$. These conditions also make $\dd\omega$ and $\delta\omega$ vanish separately. The converse and the dimension formulas follow.
\end{proof}
\begin{corollary}\label{cor:paired}
If $\xi_1,\ldots,\xi_k$ is a basis of $K$, an explicit basis of the full kernel is
\begin{equation}\label{eq:paired}
 1,\ \vol,\ \xi_1,*\xi_1,\ldots,\xi_k,*\xi_k.
\end{equation}
For $\xi=Ae^1+Be^2+Ce^3$, the paired two-form is
\begin{equation}\label{eq:pairformula}
 *\xi=Ae^{23}-Be^{13}-Ce^{12}.
\end{equation}
\end{corollary}
Complexification preserves the displayed dimensions as complex dimensions. The same reduction proof works in Riemannian signature, with the corresponding Hodge signs, but the following case-by-case results use signature $(+,+,-)$ throughout.

\section{The seven families and the classification statement}\label{sec:classification}
\subsection{Brackets, parameter restrictions, and scope}
We use the following brackets in the frame \eqref{eq:metric}. All parameters are real. The labels refer to \emph{families of metric Lie algebras}, not to seven mutually distinct abstract groups. The first four families are
\begin{align}
\g_1:\quad [e_1,e_2]&=\alpha e_1-\beta e_3,\qquad
 [e_1,e_3]=-\alpha e_1-\beta e_2,\notag\\
 [e_2,e_3]&=\beta e_1+\alpha e_2+\alpha e_3,\qquad\alpha\ne0;\label{eq:g1}\\
\g_2:\quad [e_1,e_2]&=\gamma e_2-\beta e_3,\qquad
 [e_1,e_3]=-\beta e_2-\gamma e_3,\notag\\
 [e_2,e_3]&=\alpha e_1,\qquad\gamma\ne0;\label{eq:g2}\\
\g_3:\quad [e_1,e_2]&=-\gamma e_3,\qquad
 [e_1,e_3]=-\beta e_2,\qquad [e_2,e_3]=\alpha e_1;\label{eq:g3}\\
\g_4:\quad [e_1,e_2]&=-e_2+(2\eta-\beta)e_3,\qquad
 [e_1,e_3]=-\beta e_2+e_3,\notag\\
 [e_2,e_3]&=\alpha e_1,\qquad\eta\in\{1,-1\}.\label{eq:g4}
\end{align}
No nonvanishing assumptions are imposed on $\alpha,\beta,\gamma$ in $\g_3$. The non-unimodular families are
\begin{align}
\g_5:\quad [e_1,e_2]&=0,\qquad
 [e_1,e_3]=\alpha e_1+\beta e_2,\qquad
 [e_2,e_3]=\gamma e_1+\delta e_2;\label{eq:g5}\\
\g_6:\quad [e_1,e_2]&=\alpha e_2+\beta e_3,\qquad
 [e_1,e_3]=\gamma e_2+\delta e_3,\qquad
 [e_2,e_3]=0;\label{eq:g6}\\
\g_7:\quad [e_1,e_2]&=-\alpha e_1-\beta e_2-\beta e_3,\qquad
 [e_1,e_3]=\alpha e_1+\beta e_2+\beta e_3,\notag\\
 [e_2,e_3]&=\gamma e_1+\delta e_2+\delta e_3.\label{eq:g7}
\end{align}
For all three we assume
\begin{equation}\label{eq:trace}
 t=\alpha+\delta\ne0.
\end{equation}
The customary additional restrictions are
\begin{equation}\label{eq:normalizations}
 \g_5:\ \alpha\gamma+\beta\delta=0,\qquad
 \g_6:\ \alpha\gamma-\beta\delta=0,\qquad
 \g_7:\ \alpha\gamma=0.
\end{equation}
Our formulas hold on these customary parameter sets and remain valid if any of the three restrictions in \eqref{eq:normalizations} is deleted. The Jacobi identities for \eqref{eq:g5}--\eqref{eq:g7} do not require them: each algebra has a codimension-one abelian ideal, equal in $\g_7$ to $\spn\{e_1,e_2+e_3\}$. The signs in \eqref{eq:g2} agree with \cite[\S2]{Wang}. Replacing $\gamma$ by $-\gamma$, as in another convention, does not change its kernel statements.

\subsection{Explicit notation for the answer}
For one-forms $\xi_1,\ldots,\xi_k$, define
\begin{equation}\label{eq:Bnotation}
 \cB(\xi_1,\ldots,\xi_k)
 =\spn\{1,\vol,\xi_1,*\xi_1,\ldots,\xi_k,*\xi_k\},\qquad
 \cU=\spn\{1,\vol\}.
\end{equation}
Every Hodge partner is explicitly given by \eqref{eq:pairformula}; the next two sections also write the bases without stars. In rank-one cases set
\begin{equation}\label{eq:rankone}
 (A,B)=
 \begin{cases}
 (\beta,-\alpha),&(\alpha,\beta)\ne(0,0),\\
 (\delta,-\gamma),&(\alpha,\beta)=(0,0).
 \end{cases}
\end{equation}
This definition is used only when $\alpha\delta-\beta\gamma=0$ and $t\ne0$. It always gives a nonzero pair, since in the second branch $\delta=t\ne0$.

\begin{theorem}[Classification in the seven families]\label{thm:list}
For \eqref{eq:g1}--\eqref{eq:g7}, put $J=\alpha\delta-\beta\gamma$. The full left-invariant kernel of $D$ is given in Table~\ref{tab:main}. In the $\g_3$ row, $Z\subset\{1,2,3\}$ contains index $1$ if $\alpha=0$, index $2$ if $\beta=0$, and index $3$ if $\gamma=0$. In every row the dimensions by degree are $(1,k,k,1)$, where $k$ is the number of independent displayed one-form generators.
\end{theorem}
\begin{table}[htbp]\centering\small
\renewcommand{\arraystretch}{1.17}
\begin{tabular}{@{}lll c@{}}
\toprule
Family & Parameter stratum & Full kernel & Dimension\\
\midrule
$\g_1$ & $\beta\ne0$ & $\cU$ & $2$\\
 & $\beta=0$ & $\cB(e^2-e^3)$ & $4$\\
$\g_2$ & $\alpha\ne0$ & $\cU$ & $2$\\
 & $\alpha=0$ & $\cB(e^1)$ & $4$\\
$\g_3$ & all $\alpha,\beta,\gamma$ & $\cB(e^j:j\in Z)$ & $2+2|Z|$\\
$\g_4$ & $\alpha\ne0,\ \beta\ne\eta$ & $\cU$ & $2$\\
 & $\alpha=0,\ \beta\ne\eta$ & $\cB(e^1)$ & $4$\\
 & $\alpha\ne0,\ \beta=\eta$ & $\cB(e^2+\eta e^3)$ & $4$\\
 & $\alpha=0,\ \beta=\eta$ & $\cB(e^1,e^2+\eta e^3)$ & $6$\\
$\g_5$ & $J\ne0$ & $\cU$ & $2$\\
 & $J=0$ & $\cB(Ae^1+Be^2)$ & $4$\\
$\g_6$ & $J\ne0$ & $\cU$ & $2$\\
 & $J=0$ & $\cB(Ae^2+Be^3)$ & $4$\\
$\g_7$ & $(\alpha,\gamma)\ne(0,0)$ & $\cB(e^2-e^3)$ & $4$\\
 & $\alpha=\gamma=0$ & $\cB(e^1,e^2-e^3)$ & $6$\\
\bottomrule
\end{tabular}
\caption{All parameter strata. The pair $(A,B)$ is defined in \eqref{eq:rankone}. The bracket restrictions stated in this section remain in force; no tacit genericity assumptions are made.}\label{tab:main}
\end{table}
\begin{proof}
Apply Theorem~\ref{thm:reduction} to the direct computations in Sections~\ref{sec:unimod} and \ref{sec:nonunimod}. Those calculations specify the complete operator and solve the equations in each stratum, including all exceptional parameter values.
\end{proof}

\section{Detailed computations in the unimodular families}\label{sec:unimod}
For all four families in this section $\theta=0$. Hence $\delta e^i=0$ and $\dd e^{ij}=0$. Also $D1=D\vol=0$. Consequently the formulas for $\dd$ on one-forms and $\delta$ on two-forms below specify $D$ on all eight exterior basis elements.
For a direct full-kernel calculation write
\begin{equation}\label{eq:omega8}
 \omega=f+ue^1+ve^2+we^3+pe^{12}+qe^{13}+re^{23}+h\vol.
\end{equation}
Then $D\omega=0$ is equivalent to the one-form closedness equation and the two-form coclosedness equation separately. The constants $f$ and $h$ are free. We solve both equations explicitly.

\subsection{Family $\g_1$: a null generator at $\beta=0$}
Using \eqref{eq:CE} and \eqref{eq:g1},
\begin{align}
 \dd e^1&=-\alpha e^{12}+\alpha e^{13}-\beta e^{23},\notag\\
 \dd e^2&=\beta e^{13}-\alpha e^{23},\notag\\
 \dd e^3&=\beta e^{12}-\alpha e^{23}.\label{eq:g1d}
\end{align}
The adjoint traces vanish, as can also be checked by differentiating the three two-forms. Applying $\delta_2=-*\dd*$ yields
\begin{align}
 \delta e^{12}&=-\alpha e^1-\beta e^3,\notag\\
 \delta e^{13}&=-\alpha e^1-\beta e^2,\notag\\
 \delta e^{23}&=\beta e^1+\alpha e^2-\alpha e^3.\label{eq:g1delta}
\end{align}
For example, $\delta e^{12}=-*\dd e^3=-*(\beta e^{12}-\alpha e^{23})$, which is the first line of \eqref{eq:g1delta}. Thus \eqref{eq:omega8} lies in the kernel precisely when
\begin{align}
 -\alpha u+\beta w&=0,& \alpha u+\beta v&=0,& -\beta u-\alpha(v+w)&=0,\label{eq:g1one}\\
 -\alpha(p+q)+\beta r&=0,& -\beta q+\alpha r&=0,& -\beta p-\alpha r&=0.\label{eq:g1two}
\end{align}
The coefficient matrix in \eqref{eq:g1one} is
\begin{equation}
 M_1=\begin{pmatrix}-\alpha&0&\beta\\\alpha&\beta&0\\-\beta&-\alpha&-\alpha\end{pmatrix},
 \qquad\det M_1=\beta^3.
\end{equation}
If $\beta\ne0$, all $u,v,w$ vanish. Equation \eqref{eq:g1two}, or the invertibility of the Hodge star, then gives $p=q=r=0$. Therefore $\cH_D=\cU$.

If $\beta=0$, the assumption $\alpha\ne0$ gives $u=0$, $w=-v$, $r=0$, and $q=-p$. Consequently
\begin{equation}\label{eq:g1answer}
 \omega=f+h\vol+v(e^2-e^3)+p(e^{12}-e^{13}),
\end{equation}
and an explicit basis is
\begin{equation}
 \{1,e^{123},e^2-e^3,e^{12}-e^{13}\}.
\end{equation}
The one-form $e^2-e^3$ is null. Its appearance does not require a limiting argument: it is the exact solution throughout the $\beta=0$ stratum.

\subsection{Family $\g_2$: the exceptional value $\alpha=0$}
The differential is
\begin{equation}\label{eq:g2d}
 \dd e^1=-\alpha e^{23},\qquad
 \dd e^2=-\gamma e^{12}+\beta e^{13},\qquad
 \dd e^3=\beta e^{12}+\gamma e^{13}.
\end{equation}
The codifferential on two-forms is
\begin{equation}\label{eq:g2delta}
 \delta e^{12}=-\gamma e^2-\beta e^3,\qquad
 \delta e^{13}=-\beta e^2+\gamma e^3,\qquad
 \delta e^{23}=\alpha e^1.
\end{equation}
Thus the one-form equations are
\begin{equation}\label{eq:g2one}
 \alpha u=0,\qquad -\gamma v+\beta w=0,\qquad \beta v+\gamma w=0,
\end{equation}
while the two-form equations are
\begin{equation}\label{eq:g2two}
 \alpha r=0,\qquad -\gamma p-\beta q=0,\qquad -\beta p+\gamma q=0.
\end{equation}
Since $\beta^2+\gamma^2>0$ by $\gamma\ne0$, both transverse two-by-two systems are invertible. Therefore $v=w=p=q=0$. If $\alpha\ne0$, also $u=r=0$, and $\cH_D=\cU$. If $\alpha=0$, the coefficients $u,r$ are free and
\begin{equation}\label{eq:g2answer}
 \cH_D=\spn\{1,e^{123},e^1,e^{23}\}.
\end{equation}
There is no additional branch at $\beta=0$, because $\gamma\ne0$ already guarantees invertibility of the transverse equations.

\subsection{Family $\g_3$: all eight zero/nonzero strata}
Here the differential is diagonal up to exterior ordering:
\begin{equation}\label{eq:g3d}
 \dd e^1=-\alpha e^{23},\qquad
 \dd e^2=\beta e^{13},\qquad
 \dd e^3=\gamma e^{12}.
\end{equation}
Equation \eqref{eq:delta} gives
\begin{equation}\label{eq:g3delta}
 \delta e^{12}=-\gamma e^3,\qquad
 \delta e^{13}=-\beta e^2,\qquad
 \delta e^{23}=\alpha e^1.
\end{equation}
All kernel equations therefore reduce to
\begin{equation}\label{eq:g3system}
 \alpha u=\beta v=\gamma w=0,\qquad
 \gamma p=\beta q=\alpha r=0.
\end{equation}
The constants $f,h$ are arbitrary. In particular, no condition involving a sum such as $\alpha+\beta-\gamma$ appears; these are not the spin Dirac equations.

For clarity, every parameter stratum is displayed below. Each row includes $1$ and $e^{123}$.
\begin{center}\small
\renewcommand{\arraystretch}{1.18}
\begin{tabular}{@{}lll c@{}}
\toprule
Zero parameters & Nonzero parameters & Additional basis elements & Dimension\\
\midrule
none & $\alpha,\beta,\gamma$ & none & $2$\\
$\alpha$ & $\beta,\gamma$ & $e^1,e^{23}$ & $4$\\
$\beta$ & $\alpha,\gamma$ & $e^2,e^{13}$ & $4$\\
$\gamma$ & $\alpha,\beta$ & $e^3,e^{12}$ & $4$\\
$\alpha,\beta$ & $\gamma$ & $e^1,e^2,e^{13},e^{23}$ & $6$\\
$\alpha,\gamma$ & $\beta$ & $e^1,e^3,e^{12},e^{23}$ & $6$\\
$\beta,\gamma$ & $\alpha$ & $e^2,e^3,e^{12},e^{13}$ & $6$\\
$\alpha,\beta,\gamma$ & none & $e^1,e^2,e^3,e^{12},e^{13},e^{23}$ & $8$\\
\bottomrule
\end{tabular}
\end{center}
The last row is abelian. When exactly two parameters vanish, one obtains the corresponding diagonal-bracket Heisenberg examples. Multiplication of a basis vector by $-1$ does not change its span; signs of some Hodge partners have been suppressed in this table only.

A direct check of a six-dimensional branch is useful. If $\alpha=\beta=0$ and $\gamma\ne0$, then the only nonzero actions of $D$ are $D e^3=\gamma e^{12}$ and $D e^{12}=-\gamma e^3$. The other six basis forms are killed by $D$, exactly as listed. If all three parameters vanish, $D$ is identically zero.

\subsection{Family $\g_4$: two independent degeneracy conditions}
The structure equations are
\begin{equation}\label{eq:g4d}
 \dd e^1=-\alpha e^{23},\qquad
 \dd e^2=e^{12}+\beta e^{13},\qquad
 \dd e^3=(\beta-2\eta)e^{12}-e^{13}.
\end{equation}
The corresponding two-form codifferential is
\begin{equation}\label{eq:g4delta}
 \delta e^{12}=e^2+(2\eta-\beta)e^3,\qquad
 \delta e^{13}=-\beta e^2-e^3,\qquad
 \delta e^{23}=\alpha e^1.
\end{equation}
The one-form equations are
\begin{equation}\label{eq:g4one}
 \alpha u=0,\qquad v+(\beta-2\eta)w=0,\qquad \beta v-w=0.
\end{equation}
Substituting $w=\beta v$ into the middle equation and using $\eta^2=1$ gives
\begin{equation}\label{eq:g4discriminant}
 \bigl(1+\beta(\beta-2\eta)\bigr)v=(\beta-\eta)^2v=0.
\end{equation}
Similarly, the two-form equations are
\begin{equation}\label{eq:g4two}
 \alpha r=0,\qquad p-\beta q=0,\qquad (2\eta-\beta)p-q=0.
\end{equation}
If $\beta\ne\eta$, then $v=w=p=q=0$. Thus $\alpha\ne0$ gives $\cU$, whereas $\alpha=0$ gives
\begin{equation}
 \spn\{1,e^{123},e^1,e^{23}\}.
\end{equation}
If $\beta=\eta$, then $w=\eta v$ and $p=\eta q$. The null generator and its partner are
\begin{equation}\label{eq:g4forms}
 \nu=e^2+\eta e^3,\qquad *\nu=-e^{13}-\eta e^{12}.
\end{equation}
For $\alpha\ne0$ the kernel is
\begin{equation}\label{eq:g4answer1}
 \spn\{1,e^{123},e^2+\eta e^3,e^{13}+\eta e^{12}\}.
\end{equation}
For $\alpha=0$, add $e^1$ and $e^{23}$:
\begin{equation}\label{eq:g4answer2}
 \spn\{1,e^{123},e^1,e^2+\eta e^3,e^{23},e^{13}+\eta e^{12}\}.
\end{equation}
These four possibilities exhaust the parameter set for each choice $\eta=\pm1$. In particular, setting only $\alpha=0$ does not by itself create the null generator: the additional equation $\beta=\eta$ is necessary.

\section{Detailed computations in the non-unimodular families}\label{sec:nonunimod}
We retain $t=\alpha+\delta\ne0$ and $J=\alpha\delta-\beta\gamma$. To avoid ambiguity in the calculations below, when the scalar parameter $\delta$ multiplies a form we put it in parentheses, for example $(\delta)e^1$; $\delta(e^1)$ denotes application of the codifferential.

The scalar and top-degree equations must now be included. Omitting them would incorrectly enlarge the kernel. We give $\dd$ on one-forms and two-forms, and $\delta$ on one-forms and two-forms. With $D1=D\vol=0$, these data specify the entire operator.

\subsection{Family $\g_5$: spacelike abelian ideal}
For \eqref{eq:g5},
\begin{equation}\label{eq:g5d}
 \dd e^1=-\alpha e^{13}-\gamma e^{23},\qquad
 \dd e^2=-\beta e^{13}-(\delta)e^{23},\qquad
 \dd e^3=0.
\end{equation}
The only nonzero adjoint trace is $\tr\ad_{e_3}=-t$, so
\begin{equation}\label{eq:g5mod}
 \theta=-t e^3,\qquad
 \delta(e^1)=\delta(e^2)=0,\quad \delta(e^3)=t.
\end{equation}
The last sign uses the timelike entry $g^{33}=-1$. Direct differentiation gives
\begin{equation}\label{eq:g5d2}
 \dd e^{12}=t\vol,\qquad \dd e^{13}=\dd e^{23}=0.
\end{equation}
Indeed,
\[
 \dd(e^1\wedge e^2)=
 (-\alpha e^{13}-\gamma e^{23})\wedge e^2
 -e^1\wedge(-\beta e^{13}-(\delta)e^{23})
 =(\alpha+\delta)e^{123}.
\]
Using \eqref{eq:delta},
\begin{equation}\label{eq:g5delta2}
 \delta(e^{12})=0,\qquad
 \delta(e^{13})=-(\delta)e^1+\beta e^2,\qquad
 \delta(e^{23})=\gamma e^1-\alpha e^2.
\end{equation}

For the one-form part of \eqref{eq:omega8}, the scalar equation $tw=0$ forces $w=0$. Closedness then gives
\begin{equation}\label{eq:g5one}
 \begin{pmatrix}\alpha&\beta\\\gamma&\delta\end{pmatrix}
 \begin{pmatrix}u\\v\end{pmatrix}=0.
\end{equation}
For the two-form part, the top-degree equation $tp=0$ forces $p=0$. Coclosedness gives
\begin{equation}\label{eq:g5two}
 -\delta q+\gamma r=0,\qquad \beta q-\alpha r=0.
\end{equation}
If $J\ne0$, both systems are invertible, and $\cH_D=\cU$.

Suppose $J=0$. Since $t\ne0$, the matrix in \eqref{eq:g5one} has rank exactly one. If $(\alpha,\beta)\ne(0,0)$, its kernel is generated by $(\beta,-\alpha)$ because
\begin{equation}
 \begin{pmatrix}\alpha&\beta\\\gamma&\delta\end{pmatrix}
 \begin{pmatrix}\beta\\-\alpha\end{pmatrix}
 =\begin{pmatrix}0\\-J\end{pmatrix}=0.
\end{equation}
If $(\alpha,\beta)=(0,0)$, the second row instead supplies the nonzero generator $(\delta,-\gamma)$. Thus, for the pair $(A,B)$ in \eqref{eq:rankone}, $K=\R(Ae^1+Be^2)$ and
\begin{equation}\label{eq:g5answer}
 \boxed{\ \cH_D=\spn\{1,e^{123},Ae^1+Be^2,Ae^{23}-Be^{13}\}\quad(J=0).\ }
\end{equation}
More explicitly, its basis in the first branch is
\begin{equation}\label{eq:g5branch}
 \{1,e^{123},\beta e^1-\alpha e^2,\beta e^{23}+\alpha e^{13}\},
\end{equation}
and its basis in the second branch is
\begin{equation}
 \{1,e^{123},\delta e^1-\gamma e^2,\delta e^{23}+\gamma e^{13}\}.
\end{equation}
The top and scalar equations exclude any further solutions. No step used $\alpha\gamma+\beta\delta=0$. Under that restriction, the two rows of the matrix are orthogonal in a Euclidean plane, so singularity forces one of them to vanish; our two-branch answer includes both possibilities.

\subsection{Family $\g_6$: Lorentzian abelian ideal}
For \eqref{eq:g6},
\begin{equation}\label{eq:g6d}
 \dd e^1=0,\qquad
 \dd e^2=-\alpha e^{12}-\gamma e^{13},\qquad
 \dd e^3=-\beta e^{12}-(\delta)e^{13}.
\end{equation}
The nonzero adjoint trace is $\tr\ad_{e_1}=t$, so
\begin{equation}\label{eq:g6mod}
 \theta=t e^1,\qquad \delta(e^1)=t,\quad \delta(e^2)=\delta(e^3)=0.
\end{equation}
The two-form differential is
\begin{equation}\label{eq:g6d2}
 \dd e^{12}=\dd e^{13}=0,\qquad \dd e^{23}=-t\vol.
\end{equation}
For example, the last identity follows from
$\dd e^2\wedge e^3-e^2\wedge\dd e^3
=-\alpha e^{123}-(\delta)e^{123}$.
The codifferential on two-forms is
\begin{equation}\label{eq:g6delta2}
 \delta(e^{12})=(\delta)e^2+\beta e^3,\qquad
 \delta(e^{13})=\gamma e^2+\alpha e^3,\qquad \delta(e^{23})=0.
\end{equation}
The one-form equations force $u=0$, followed by
\begin{equation}\label{eq:g6one}
 \alpha v+\beta w=0,\qquad \gamma v+\delta w=0.
\end{equation}
The top-degree equation forces $r=0$, and the two-form equations become
\begin{equation}\label{eq:g6two}
 \delta p+\gamma q=0,\qquad \beta p+\alpha q=0.
\end{equation}
Thus $J\ne0$ gives $\cH_D=\cU$.

If $J=0$, the pair $(A,B)$ in \eqref{eq:rankone} gives $K=\R(Ae^2+Be^3)$. Its Hodge partner is $-Ae^{13}-Be^{12}$. An explicit basis is therefore
\begin{equation}\label{eq:g6answer}
 \boxed{\ \{1,e^{123},Ae^2+Be^3,Ae^{13}+Be^{12}\}\quad(J=0).\ }
\end{equation}
In the first branch it becomes
\begin{equation}
 \{1,e^{123},\beta e^2-\alpha e^3,\beta e^{13}-\alpha e^{12}\},
\end{equation}
and in the second branch it is
\begin{equation}
 \{1,e^{123},\delta e^2-\gamma e^3,\delta e^{13}-\gamma e^{12}\}.
\end{equation}
The one-form generator has squared norm $A^2-B^2$ and can be spacelike, timelike, or null. Its causal type does not change the dimension, since the span is still one-dimensional. Neither the equations nor their solution use $\alpha\gamma-\beta\delta=0$.

To check a null singular branch explicitly, take $\alpha=\beta=\gamma=\delta=1$. Then $t=2$, $J=0$, and the customary restriction holds. The two kernel generators are $e^2-e^3$ and $e^{13}-e^{12}$. By contrast, if $\alpha=\delta=1$, $\beta=0$, $\gamma=1$, then $J=1$; even though some coefficients vanish, both middle-degree kernels are zero.

\subsection{Family $\g_7$: degenerate abelian ideal}
Introduce the null form and its partner
\begin{equation}\label{eq:NandT}
 N=e^2-e^3,\qquad T=e^{12}-e^{13}=*N.
\end{equation}
The structure equations are
\begin{equation}\label{eq:g7d}
 \dd e^1=\alpha T-\gamma e^{23},\qquad
 \dd e^2=\beta T-(\delta)e^{23},\qquad
 \dd e^3=\beta T-(\delta)e^{23}.
\end{equation}
In particular $\dd N=0$, without any parameter restriction. The modular character and one-form codifferential are
\begin{equation}\label{eq:g7mod}
 \theta=tN,\qquad
 \delta(e^1)=0,\quad \delta(e^2)=t,\quad \delta(e^3)=t.
\end{equation}
Both last signs are positive because $\theta_3=-t$ and $g^{33}=-1$. Differentiation of the two-forms gives
\begin{equation}\label{eq:g7d2}
 \dd e^{12}=t\vol,\qquad \dd e^{13}=t\vol,\qquad \dd e^{23}=0.
\end{equation}
The codifferential is
\begin{align}
 \delta(e^{12})&=-(\delta)e^1+\beta(e^2-e^3),\notag\\
 \delta(e^{13})&=-(\delta)e^1+\beta(e^2-e^3),\notag\\
 \delta(e^{23})&=\gamma e^1-\alpha(e^2-e^3).\label{eq:g7delta2}
\end{align}
In particular $\delta T=0$, and $DN=DT=0$.

For the one-form component, the scalar and two-form equations are
\begin{align}
 t(v+w)&=0,\notag\\
 \alpha u+\beta(v+w)&=0,\qquad
 \gamma u+\delta(v+w)=0.\label{eq:g7one}
\end{align}
Since $t\ne0$, they reduce to
\begin{equation}\label{eq:g7onereduced}
 w=-v,\qquad \alpha u=\gamma u=0.
\end{equation}
For the two-form component, $t(p+q)=0$, and the one-form equations are
\begin{equation}\label{eq:g7two}
 -\delta(p+q)+\gamma r=0,\qquad
 \beta(p+q)-\alpha r=0.
\end{equation}
Consequently $q=-p$ and $\gamma r=\alpha r=0$.

If $(\alpha,\gamma)\ne(0,0)$, then $u=r=0$ and
\begin{equation}\label{eq:g7answer1}
 \boxed{\ \cH_D=\spn\{1,e^{123},e^2-e^3,e^{12}-e^{13}\}.\ }
\end{equation}
If $\alpha=\gamma=0$, then $\delta=t\ne0$, while $u,r$ are free. Hence
\begin{equation}\label{eq:g7answer2}
 \boxed{\ \cH_D=\spn\{1,e^{123},e^1,e^2-e^3,e^{23},e^{12}-e^{13}\}.\ }
\end{equation}
There is no further exceptional condition on $\beta$. No step used $\alpha\gamma=0$; in particular, when $\alpha\gamma\ne0$, the unrestricted family still has the four-dimensional kernel \eqref{eq:g7answer1}. This is why the computation remains valid beyond the customary normalization.

\section{Invariant interpretation and coverage beyond normalization}\label{sec:geometry}
\subsection{The modular direction and dimension jumps}
Let $b=\dim\Ann([\g,\g])$ and let $r_\theta$ be the rank of the functional $\langle\theta,\cdot\rangle$ restricted to $W=\Ann([\g,\g])$. For a non-unimodular algebra, $r_\theta$ is zero or one. Theorem~\ref{thm:reduction} gives
\begin{equation}\label{eq:dimrank}
 \dim\cH_D=2+2(b-r_\theta).
\end{equation}
Here $b$ is also the first Lie-algebra Betti number, since the differential of invariant zero-forms is zero. It is not being identified with the first Betti number of an arbitrary underlying manifold or quotient.

\begin{proposition}\label{prop:modular-rank}
The functional $\langle\theta,\cdot\rangle|_W$ vanishes if and only if
\begin{equation}\label{eq:thetaderived}
 \theta^\sharp\in[\g,\g].
\end{equation}
When $\theta\ne0$, this condition implies that $\theta$ is null. If $b=1$, nullity is sufficient as well. If $b=2$, nullity need not suffice.
\end{proposition}
\begin{proof}
For $\xi\in W$, $\langle\theta,\xi\rangle=\xi(\theta^\sharp)$. Vanishing for all $\xi\in\Ann([\g,\g])$ is precisely \eqref{eq:thetaderived}. Since $\theta\in W$, it implies $\langle\theta,\theta\rangle=0$. If $W=\R\theta$, this one condition suffices. When $W$ is two-dimensional, it tests only one vector in a two-dimensional domain.
\end{proof}
For $\g_5$, the modular character $-te^3$ is timelike; for $\g_6$, the modular character $te^1$ is spacelike. Thus $r_\theta=1$ in both families. The derived algebra has dimension two for $J\ne0$ and one for $J=0$, yielding total kernel dimensions two and four.

For $\g_7$, $\theta=tN$ is always null, which explains the persistent four-dimensional kernel. The further increase to dimension six occurs precisely when $\alpha=\gamma=0$. In this case the brackets all take values in $\R(e_2+e_3)$ and
$\theta^\sharp=t(e_2+e_3)\in[\g,\g]$. If $J=0$ but $(\alpha,\gamma)\ne(0,0)$, the derived algebra is still one-dimensional but does not contain that modular vector, so the dimension remains four. Thus a singular bracket matrix alone does not determine the kernel in the degenerate-ideal case.

For a concrete comparison within $\g_7$, choose $\beta=\gamma=0$ and $\delta=1$. If $\alpha=0$, the kernel is six-dimensional. If $\alpha=1$, then $J=1$ and it is four-dimensional. A singular four-dimensional example is $\alpha=1$, $\beta=\delta=0$ and $\gamma=1$ in the unrestricted family: $t=1$, $J=0$, yet the condition $(\alpha,\gamma)\ne(0,0)$ still gives only four dimensions.

\subsection{A normalization-independent covering argument}
The non-unimodular classification can be used without any extra row-orthogonality assumption.
\begin{lemma}\label{lem:almostabelian}
Every three-dimensional non-unimodular real Lie algebra has a two-dimensional abelian ideal $\mathfrak u=\ker\theta$.
\end{lemma}
\begin{proof}
Equation \eqref{eq:theta-closed} makes $\mathfrak u$ a two-dimensional ideal. For $X\in\mathfrak u$, the induced action of $\ad_X$ on $\g/\mathfrak u$ is zero, since $[X,\g]\subset\mathfrak u$. Hence $\tr(\ad_X|_{\mathfrak u})=\theta(X)=0$. A nonabelian two-dimensional algebra has a basis $(a,b)$ with $[a,b]=b$, and its adjoint $\ad_a$ has trace one. Thus $\mathfrak u$ is abelian.
\end{proof}
If the metric on $\mathfrak u$ is positive definite, choose a spacelike orthonormal basis $e_1,e_2$ in it and a timelike unit normal $e_3$. The brackets have the unrestricted form $\g_5$. If $\mathfrak u$ is Lorentzian, choose $e_2,e_3$ in it and a spacelike unit normal $e_1$; the brackets have the unrestricted form $\g_6$.

If $\mathfrak u$ is degenerate, choose a unit spacelike $x\in\mathfrak u$, a null radical generator $y\in\mathfrak u$, and a transverse null vector $z$ with $g(y,z)=1$ and $g(x,z)=0$. Write
\begin{equation}\label{eq:nullbrackets}
 [z,x]=a x+c y,\qquad [z,y]=b x+d y,\qquad [x,y]=0,\qquad a+d\ne0.
\end{equation}
Set
\begin{equation}\label{eq:nullframe}
 e_1=x,\qquad e_2=\tfrac12 y+z,\qquad e_3=\tfrac12 y-z.
\end{equation}
This is a pseudo-orthonormal frame of signature $(+,+,-)$, with $e_2+e_3=y$. Its brackets are
\begin{align}
 [e_1,e_2]&=-a e_1-c(e_2+e_3),\notag\\
 [e_1,e_3]&=a e_1+c(e_2+e_3),\notag\\
 [e_2,e_3]&=b e_1+d(e_2+e_3).\label{eq:nullactual}
\end{align}
The last sign follows explicitly from
$[\tfrac12y+z,\tfrac12y-z]=[z,y]$.
Thus \eqref{eq:nullactual} is $\g_7$ with $\alpha=a$, $\beta=c$, $\gamma=b$, $\delta=d$. The trace condition remains $t=a+d\ne0$.

\begin{corollary}\label{cor:coverage}
Every three-dimensional non-unimodular Lorentzian metric Lie algebra can be represented by one of the unrestricted bracket families $\g_5,\g_6,\g_7$. The kernel formulas in Section~\ref{sec:nonunimod} cover all such representatives without the extra restrictions \eqref{eq:normalizations}.
\end{corollary}
The four unimodular normal-form types are the classical Lorentzian self-adjoint bracket-operator types; see \cite{Rahmani,BatatOnda}. Together with the above covering argument, this yields a kernel computation that is not affected by incompleteness of some customary non-unimodular normalizations. It is not necessary to assert uniqueness of these representatives, nor to quotient the displayed bases by automorphisms.

For example, in unrestricted $\g_6$ take $\alpha=\delta=1$, $\beta=1$, $\gamma=-1$. The normalization expression $\alpha\gamma-\beta\delta$ equals $-2$, while $J=2$. Its kernel is nevertheless immediately $\cU$. In unrestricted $\g_7$, the choice $\alpha=\gamma=\delta=1$, $\beta=0$ has $\alpha\gamma\ne0$ and the four-dimensional basis \eqref{eq:g7answer1}. These inputs illustrate how the formulas can be applied before a normal-form restriction is imposed.

\subsection{Comparison with the spin Dirac problem}
The distinction from \cite{GilRusso} is already visible in the universal pair $1,\vol$. A spinor zero mode is not a differential form, and the exterior module is not the irreducible spin module. Scalar-curvature restrictions and special metric relations obtained from a spin Dirac matrix therefore cannot be transferred to \eqref{eq:Dintro}.

For any three-dimensional simple Lie algebra, $[\g,\g]=\g$, so $W=0$ and $\cH_D=\cU$ for every Lorentzian metric. The existence of this two-dimensional form kernel says nothing about the existence of invariant spinor zero modes. Nilpotent examples in $\g_3$, by contrast, have six-dimensional exterior kernels. The shared feature is a reduction to structure constants, not equality of the two kernel problems.

Our calculation also avoids a global assumption: no spin structure or simple connectivity is needed to define invariant differential forms. On a noncompact group these forms need not be square-integrable. If they descend to a quotient, their algebraic equations remain valid, but the complete kernel on that quotient may contain non-invariant forms. No such larger kernel is classified here.

\section{The first-order kernel and the kernel of its square}\label{sec:square}
Although $\ker D\subseteq\ker\Delta$ for $\Delta=D^2$, equality is not automatic. The next example lies entirely within one of the seven families.
\begin{proposition}\label{prop:jordan}
In $\g_1$ with $\alpha=1$, $\beta=0$,
\begin{equation}\label{eq:jordan}
 \dim\ker D=4,\qquad \dim\ker D^2=6.
\end{equation}
\end{proposition}
\begin{proof}
Put $N=e^2-e^3$ and $T=e^{12}-e^{13}$. Equations \eqref{eq:g1d}--\eqref{eq:g1delta} give
\begin{align}
 D e^1&=-T,&D e^2&=-e^{23},&D e^3&=-e^{23},\notag\\
 D e^{12}&=-e^1,&D e^{13}&=-e^1,&D e^{23}&=N,\label{eq:jordan-action}
\end{align}
with $DN=DT=0$. Squaring,
\begin{align}
 D^2e^1&=0,&D^2e^2&=-N,&D^2e^3&=-N,\notag\\
 D^2e^{12}&=T,&D^2e^{13}&=T,&D^2e^{23}&=0.\label{eq:jordan-square}
\end{align}
The images $N$ and $T$ are nonzero and have distinct degrees, so $D^2$ has rank two and its kernel has dimension six. An explicit basis is
\begin{equation}
 \{1,e^{123},N,T,e^1,e^{23}\}.
\end{equation}
The last two forms are not $D$-harmonic. Theorem~\ref{thm:list}, or direct inspection of \eqref{eq:jordan-action}, gives $\dim\ker D=4$.
\end{proof}

There is an even more degenerate non-unimodular example. In $\g_7$, let $\alpha=\beta=\gamma=0$ and $\delta=1$. Then
\begin{align}
 D e^1&=0,&D e^2&=1-e^{23},&D e^3&=1-e^{23},\notag\\
 D e^{12}&=-e^1+\vol,&D e^{13}&=-e^1+\vol,&D e^{23}&=0.
\label{eq:nilD}
\end{align}
Together with $D1=D\vol=0$, these equations give $D^2=0$ on the entire invariant exterior algebra, although $D\ne0$. Thus $\dim\ker D=6$ but $\dim\ker D^2=8$. This example satisfies the customary condition $\alpha\gamma=0$ and requires no extension of the standard family.

These examples also explain the terminology adopted here. In Lorentzian geometry a form satisfying $\Delta\omega=0$ is often called harmonic, whereas the first-order equation may be called closed and coclosed, or a Hodge--de Rham zero mode. Our classification always concerns the first-order equation; Theorem~\ref{thm:reduction} proves that it is precisely the closed-and-coclosed condition on invariant forms in dimension three. It does not identify it with the wave-type equation $\Delta\omega=0$.

\section{Concluding remarks}
After fixing the orientation and geometric codifferential, the complete invariant kernel problem reduces to $\dd\xi=0$ and $\langle\theta,\xi\rangle=0$ for one-forms. Table~\ref{tab:main} and Sections~\ref{sec:unimod}--\ref{sec:nonunimod} provide concrete bases on every parameter stratum of the seven bracket families.

The unrestricted non-unimodular formulas separate the kernel computation from the choice of metric normalization. They show that dimension jumps depend on the derived algebra and modular direction, rather than on a particular expression for the Levi-Civita connection. The paired degree-one and degree-two generators make this dependence explicit.

Further questions concern the complete smooth or $L^2$ kernel on a group or a quotient, and the larger kernel of the Hodge Laplacian. Those are different analytic problems. In higher dimensions inhomogeneous components may also cancel in ways excluded by \eqref{eq:Dinhom}, so a direct generalization of the three-dimensional splitting requires additional arguments.

\appendix
\section{Full matrix data and independent verification}\label{sec:appendix}
\subsection{The universal eight-dimensional matrix}
Use the ordered exterior basis
\begin{equation}
 (1,e^1,e^2,e^3,e^{12},e^{13},e^{23},\vol).
\end{equation}
Write $\dd e^j=\sum_{I\in\{12,13,23\}}M_{Ij}e^I$. If $\theta=\theta_1e^1+\theta_2e^2+\theta_3e^3$, set
\begin{equation}\label{eq:ml}
 m=(\theta_1,\theta_2,-\theta_3),\quad
 \ell=(-\theta_3,\theta_2,-\theta_1),\quad
 S=\begin{pmatrix}0&0&-1\\0&1&0\\1&0&0\end{pmatrix}.
\end{equation}
Here $S$ is the Hodge star from two-forms to one-forms in the specified orders; the inverse-direction star has matrix $-S^{-1}$. The two-form codifferential is
\begin{equation}\label{eq:Cmatrix}
 C=-SMS.
\end{equation}
With block sizes $(1,3,3,1)$, the complete operator and its square are
\begin{equation}\label{eq:8matrix}
 [D]=\begin{pmatrix}
 0&m&0&0\\
 0&0&C&0\\
 0&M&0&0\\
 0&0&\ell&0
 \end{pmatrix},\qquad
 [D^2]=\operatorname{diag}(0,CM,MC,0).
\end{equation}
The potentially off-diagonal products $\ell M$ and $mC$ vanish by $\dd^2=\delta^2=0$.

The seven matrices are
\begin{align}
 M_1&=\begin{pmatrix}-\alpha&0&\beta\\\alpha&\beta&0\\-\beta&-\alpha&-\alpha\end{pmatrix},&
 M_2&=\begin{pmatrix}0&-\gamma&\beta\\0&\beta&\gamma\\-\alpha&0&0\end{pmatrix},\label{eq:M12}\\
 M_3&=\begin{pmatrix}0&0&\gamma\\0&\beta&0\\-\alpha&0&0\end{pmatrix},&
 M_4&=\begin{pmatrix}0&1&\beta-2\eta\\0&\beta&-1\\-\alpha&0&0\end{pmatrix},\label{eq:M34}\\
 M_5&=\begin{pmatrix}0&0&0\\-\alpha&-\beta&0\\-\gamma&-\delta&0\end{pmatrix},&
 M_6&=\begin{pmatrix}0&-\alpha&-\beta\\0&-\gamma&-\delta\\0&0&0\end{pmatrix},\label{eq:M56}\\
 M_7&=\begin{pmatrix}\alpha&\beta&\beta\\-\alpha&-\beta&-\beta\\-\gamma&-\delta&-\delta\end{pmatrix}.&&\label{eq:M7}
\end{align}
For $i=1,2,3,4$, one has $m=\ell=0$. For the last three,
\begin{equation}\label{eq:mlcases}
\begin{array}{c|cc}
 &m&\ell\\\hline
\g_5&(0,0,t)&(t,0,0)\\
\g_6&(t,0,0)&(0,0,-t)\\
\g_7&(0,t,t)&(t,t,0).
\end{array}
\end{equation}
Equations \eqref{eq:8matrix}--\eqref{eq:mlcases} specify all seven complete operators at both generic and singular parameter values.

\subsection{The finite-dimensional adjoint correction}
In these bases the Gram matrices on one-forms and two-forms are
\begin{equation}
 Q_1=\operatorname{diag}(1,1,-1),\qquad
 Q_2=\operatorname{diag}(1,-1,-1).
\end{equation}
The finite-dimensional adjoint of $\dd:\Lambda^1\to\Lambda^2$ is $Q_1^{-1}M^{\mathsf T}Q_2$. Contraction on two-forms by $\theta^\sharp$ has matrix
\begin{equation}
 I_\theta=\begin{pmatrix}
 -\theta_2&\theta_3&0\\
 \theta_1&0&\theta_3\\
 0&\theta_1&\theta_2
 \end{pmatrix}.
\end{equation}
Writing $\dd e^j=M_{1j}e^{12}+M_{2j}e^{13}+M_{3j}e^{23}$, the modular components are
\begin{equation}
 \theta_1=-M_{12}-M_{23},\qquad
 \theta_2=M_{11}-M_{33},\qquad
 \theta_3=M_{21}+M_{32}.
\end{equation}
Substitution now proves entry by entry that
\begin{equation}
 -SMS=Q_1^{-1}M^{\mathsf T}Q_2+I_\theta.
\end{equation}
On one-forms the finite-dimensional adjoint is zero and the contraction is the row $m$. On top forms the finite-dimensional adjoint of $\ell$ is $-\ell^{\mathsf T}$ multiplied by $Q_2^{-1}$; its sum with contraction by $\theta^\sharp$ is zero. Thus $\delta=\dd^\dagger+\iota_{\theta^\sharp}$ in every degree, as asserted in Remark~\ref{rem:adjoint}. This calculation exhibits the correction without a global integration argument.

\subsection{Rank verification and parameter branches}
For independent exact verification, construct the Chevalley--Eilenberg matrices directly from each bracket, then construct $\delta$ using the explicit star matrices \eqref{eq:star0}--\eqref{eq:star2}. The result must agree with \eqref{eq:8matrix}. In particular
\begin{equation}
 \ell_iM_i=0,\qquad m_iC_i=0
\end{equation}
are polynomial identities in every displayed family. Each proposed kernel vector is checked by multiplication by this full matrix under its stated parameter constraints.

Completeness requires parameter stratification rather than a generic rank computation. The reduced one-form system is
\begin{equation}
 \begin{pmatrix}M_i\\m_i\end{pmatrix}
 \begin{pmatrix}u\\v\\w\end{pmatrix}=0.
\end{equation}
For $\g_1$ the determinant branch is $\beta^3$; for $\g_2$ the transverse determinant has magnitude $\beta^2+\gamma^2$; for $\g_4$ it has magnitude $(\beta-\eta)^2$; for $\g_5,\g_6$ the remaining determinant is $J$; and $\g_7$ reduces to \eqref{eq:g7onereduced}. For $\g_3$, all three conditions are independent scalar equations \eqref{eq:g3system}. These calculations give the rank on each complete parameter stratum.

In a computational implementation, a useful final check is
\[
 \dim\ker[D]=2+2\dim\ker\begin{pmatrix}M_i\\m_i\end{pmatrix}.
\]
This tests the inhomogeneous kernel, not merely the individual-degree equations. The squared matrix must separately reproduce the two strict kernel inclusions in Section~\ref{sec:square}.

\end{document}